\documentstyle{article}
\begin{document}
\begin{center}
\Large{\bf  On geometrical properties of the spaces defined \\[2mm] by the Pfaff
equations
 }\vspace{4mm}\normalsize
\end{center}
 \begin{center}
\Large{\bf Valery Dryuma}\vspace{4mm}\normalsize
\end{center}
\begin{center}
{\bf Institute of Mathematics and Informatics AS Moldova,
Kishinev}\vspace{4mm}\normalsize
\end{center}
\begin{center}
{\bf  E-mail: valery@dryuma.com;\quad cainar@mail.md}\vspace{4mm}\normalsize
\end{center}
\begin{center}
{\bf  Abstract}\vspace{4mm}\normalsize
\end{center}

{Geometrical properties of holonomic and non holonomic varieties defined by the
Pfaff equations connected with the first order system of equations are studied. The
Riemann extensions of affine connected spaces for investigations of geodesics and
asymptotic lines are used}

\section{Introduction}

   There is the connection between of the Pfaff
   equation

\begin{equation} \label{dryuma:eq0}
P(x,y,z)\frac{d x }{ds}+Q(x,y,z) \frac{d y }{ds}+R(x,y,z)\frac{d z
}{ds}=0
\end{equation}

and the  first order system differential equations in form

\begin{equation} \label{dryuma:eq1}
\frac{d x }{P(x,y,z)}= \frac{d y }{Q(x,y,z)}= \frac{d z
}{R(x,y,z)}.
\end{equation}

    For example the equation (\ref{dryuma:eq0}) is exactly integrable at the
    conditions
$$ \frac{\partial P(x,y,z)}{\partial y}- \frac{\partial
Q(x,y,z)}{\partial x}=0,\quad \frac{\partial Q(x,y,z)}{\partial
z}-\frac{\partial R(x,y,z)}{\partial y}=0, $$

$$
\frac{\partial R(x,y,z)}{\partial x}-\frac{\partial P(x,y,z)}{\partial z}=0
$$
and its general integral determines the family of the surfaces in
$R^3$ space
\[
V(x,y,z)=constant,
\]
which are orthogonal to the lines of the vector field
\begin{equation} \label{dryuma:eq3}
\vec N=(P(x,y,z),Q(x,y,z),R(x,y,z)).
\end{equation}

    In more general case
    \[
P(x,y,z)\left(\frac{\partial Q(x,y,z)}{\partial z}- \frac{\partial
R(x,y,z)}{\partial y}\right)+Q(x,y,z)\left(\frac{\partial
R(x,y,z)}{\partial x}-\frac{\partial P(x,y,z)}{\partial
z}\right)+\]\[+ R(x,y,z)\left(\frac{\partial R(x,y,z)}{\partial
y}-\frac{\partial Q(x,y,z)}{\partial x}\right)=0
\]
the equation (\ref{dryuma:eq0})
\[
\mu\left(P(x,y,z)d x +Q(x,y,z) d y +R(x,y,z)d z\right)=d U(x,y,z)
\]
is also integrable by means of integrating multiplier $\mu$ and determines the
family of the surfaces $U(x,y,z)=const$ passing through the each point of the space
and orthogonal to the vector field (\ref{dryuma:eq3}).

   The vector field (\ref{dryuma:eq3}) in the space $R^3$  with conditions
\[
\left(\vec N, rot \vec N\right)=0
\]
  is called holonomic.

   In general case the Pfaff equation (\ref{dryuma:eq0}) is
   not integrable and them corresponds  the system of the integral curves
   (Pfaff variety) passing through  each point $(x,y,z)$ with tangent lines lying on the plane
\begin{equation} \label{dryuma:eq4}
P(x,y,z)(X-x)+Q(x,y,z)(Y-y)+R(x,y,z)(Z-z)=0.
\end{equation}

  The  set of the planes and points  (\ref{dryuma:eq4}) defined by the equation (\ref{dryuma:eq0}) in general case
   forms two-dimensional non holonomic variety $M^2$ which is generalization of the surface.

   For the variety $M^2$  may be extended many of the results of classical differential
   geometry of surfaces. For example the notion of the asymptotic lines, curvature lines and geodesic has
   analog for the variety $M^2$  \cite{dryuma:sintsov, dryuma4:sluchaev, dryuma5:aminov}.

   In fact the solutions of the system
\[
Pdx+Qdy+Rdz=0, \quad dP dx+dQ dy+dR dz=0
\]
or
\[
P\dot x+\dot y+R \dot z=0,
\]
\[ P_x(\dot x)^2+Q_y (\dot
y)^2+R_z(\dot z)^2+(P_y+Q_x)\dot x \dot y +(P_z+R_x)\dot x \dot
z+(Q_z+R_y)\dot y \dot z=0
\]

\noindent  give us the  curve lines of the variety $M^2$ which are the analog of the of asymptotic lines on the
   holonomic surface.

   The notion of the curvature lines also can be generalized on the variety $M^2$.

   They may be of the two kinds and one of them is defined by the solutions of the
system of equations
\[
Pdx+Qdy+Rdz=0,
\]
$$ \left|
\begin{array}{ccc}
 2P_x dx+(Q_x+P_y) dy
+(P_x+R_x) dz & P & dx \\
 (Q_x+P_y) dx + 2Q_y dy+(Q_z+R_y) dz & Q & dy \\
  (P_z+R_x) dx +(Q_z+R_y) dy +2R_z dz & R & dz
\end{array}
\right| =0 .
$$

    The notion of the geodesics  also can be extended on the variety $V^2$ (\ref{dryuma:eq0}) and they may be of
     two types.

  The first type  is determined from the condition
$$
\left [\begin {array}{ccc} P(x,y,z)&Q(x,y,z)&R(x,y,z)
\\\noalign{\medskip}\displaystyle{\frac {d}{ds}}x(s)&\displaystyle{\frac {d}{ds}}y(s)&\displaystyle{\frac {d}
{ds}}z(s)\\\noalign{\medskip}\displaystyle{\frac {d^{2}}{d{s}^{2}}}x(s)&\displaystyle{\frac
{d^{ 2}}{d{s}^{2}}}y(s)&\displaystyle{\frac {d^{2}}{d{s}^{2}}}z(s)\end
{array}\right ] =0,
$$

\noindent or
\[\left (P(x,y,z){\frac {d}{ds}}y(s)-\left ({\frac
{d}{ds}}x(s)\right )Q (x,y,z)\right ){\frac
{d^{2}}{d{s}^{2}}}z(s)+\]\[+\left (-P(x,y,z){\frac {d
}{ds}}z(s)+\left ({\frac {d}{ds}}x(s)\right )R(x,y,z)\right
){\frac {d ^{2}}{d{s}^{2}}}y(s)+\]\[+\left ({\frac
{d^{2}}{d{s}^{2}}}x(s)\right ) \left (Q(x,y,z){\frac
{d}{ds}}z(s)-R(x,y,z){\frac {d}{ds}}y(s)\right ) =0. \]
\medskip
 This relation is equivalent the system of equations

\begin{equation} \label{dryuma:eq5}
\frac{d^2 x}{ds^2}+\frac{P}{P^2+Q^2+R^2}\left(\frac{d x}{ds}\frac{
d P}{ds}+\frac{d y}{ds}\frac{ d Q}{ds}+\frac{d z}{ds}\frac{ d
R}{ds}\right)=0, \end{equation}
\medskip
\begin{equation} \label{dryuma:eq6}
\frac{d^2 y}{ds^2}+\frac{Q}{P^2+Q^2+R^2}\left(\frac{d x}{ds}\frac{
d P}{ds}+\frac{d y}{ds}\frac{ d Q}{ds}+\frac{d z}{ds}\frac{ d
R}{ds}\right)=0, \end{equation}
\medskip
\begin{equation} \label{dryuma:eq7}
 \frac{d^2 z}{ds^2}+\frac{R}{P^2+Q^2+R^2}\left(\frac{d
x}{ds}\frac{ d P}{ds}+\frac{d y}{ds}\frac{ d Q}{ds}+\frac{d
z}{ds}\frac{ d R}{ds}\right)=0.
 \end{equation}
\medskip

     Remark that after substitution of the corresponding expressions for the second derivatives
     on coordinates from the  (\ref{dryuma:eq5})--(\ref{dryuma:eq7}) into the relations
$$
\frac{dP(x,y,z)}{ds}\frac{d x }{ds}\\+\\\frac{dQ(x,y,z)}{ds}
\frac{d y }{ds}\\+\\\frac{dR(x,y,z)}{ds}\frac{d z
}{ds}+P(x,y,z)\frac{d^2 x }{ds^2}\\+
$$
$$
+\\Q(x,y,z) \frac{d^2 y
}{ds^2}\\+\\R(x,y,z)\frac{d^2 z }{ds^2}=0
$$

\noindent one get the identity.

   The definition of the second type of geodesic in the system of integral curves of
    the Pfaff equation (\ref{dryuma:eq0}) is more complicated and does not be  used
     hereinafter.

\section{The  Riemann extension of the affine connected space}.

    The formulas (\ref{dryuma:eq5})--(\ref{dryuma:eq7}) can be
    rewritten in form of geodesic of the $R^3$-space equipped with symmetrical  affine connection
    $ \Pi^i_{jk}(x^l)=\Pi^i_{kj}(x^l)$
\[
\frac{d^2 x^i}{ds^2}+\Pi^i_{jk}\frac{d x^j}{ds}\frac{d x^k}{ds}=0.
\]

   In our case we get the system of equations
$$
\frac{d^2 x}{ds^2}+\frac{P}{\Delta}\left(P_x \left(\frac{d
x}{ds}\right)^2\!+\!(P_y+Q_x) \frac{dx}{ds}\frac{d y}{ds}\!+\!(P_z+R_x)\frac{d
x}{ds} \frac{d z}{ds}\!+\!Q_y \left(\frac{ dy}{ds}\right)^2\!+\right.
$$
$$
\left.+\!(Q_z+R_y)\frac{d y}{ds}\frac{d z}{ds}\!+\!R_z \left(\frac{ d z}{ds}\right)^2\right)=0,
$$

$$
\frac{d^2 y}{ds^2}+\frac{Q}{\Delta}\left(P_x \left(\frac{d
x}{ds}\right)^2\!+\!(P_y+Q_x) \frac{dx}{ds}\frac{d y}{ds}\!+\!(P_z+R_x)\frac{d
x}{ds} \frac{d z}{ds}\!+\!Q_y \left(\frac{ dy}{ds}\right)^2\!+\right.
$$
$$
\left.+\!(Q_z+R_y)\frac{d
y}{ds}\frac{d z}{ds}\!+\!R_z \left(\frac{ d z}{ds}\right)^2\right)=0,
$$

$$
\frac{d^2 z}{ds^2}+\frac{R}{\Delta}\left(P_x \left(\frac{d
x}{ds}\right)^2\!+\!(P_y+Q_x) \frac{dx}{ds}\frac{d y}{ds}\!+\!(P_z+R_x)\frac{d
x}{ds} \frac{d z}{ds}\!+\!Q_y \left(\frac{ dy}{ds}\right)^2\!+\right.
$$
$$
\left.+\!(Q_z+R_y)\frac{d
y}{ds}\frac{d z}{ds}\!+\!R_z \left(\frac{ d z}{ds}\right)^2\right)=0,
$$
where \[\Delta=P^2+Q^2+R^2,\] from which we get the expressions
for the coefficients of affine connection. They are
$$
\Pi^1_{11}=\frac{PP_x}{\Delta},\qquad
\Pi^1_{22}=\frac{PQ_y}{\Delta},\qquad
\Pi^1_{33}=\frac{PR_z}{\Delta},
$$

$$
\Pi^1_{12}=\frac{P(P_y+Q_x)}{2\Delta},\qquad
\Pi^1_{13}=\frac{P(P_z+R_x)}{2\Delta},\qquad
\Pi^1_{23}=\frac{P(Q_z+R_y)}{2\Delta},
$$

$$
\Pi^2_{11}=\frac{QP_x}{\Delta},\qquad
\Pi^2_{22}=\frac{QQ_y}{\Delta},\qquad
\Pi^2_{33}=\frac{QR_z}{\Delta},
$$

$$
\Pi^2_{12}=\frac{Q(P_y+Q_x)}{2\Delta},\qquad
\Pi^2_{13}=\frac{Q(P_z+R_x)}{2\Delta},\qquad
\Pi^2_{23}=\frac{Q(Q_z+R_y)}{2\Delta},
$$

$$
\Pi^3_{11}=\frac{RP_x}{\Delta},\qquad
\Pi^3_{22}=\frac{RQ_y}{\Delta},\qquad
\Pi^3_{33}=\frac{RR_z}{\Delta},
$$

$$
\Pi^3_{12}=\frac{R(P_y+Q_x)}{2\Delta},\qquad
\Pi^3_{13}=\frac{R(P_z+R_x)}{2\Delta},\qquad
\Pi^3_{23}=\frac{R(Q_z+R_y)}{2\Delta}.
$$
\medskip

     So with any equation (\ref{dryuma:eq0})
     can be associated 3-dimensional affine connected space and its properties
     will be dependent from the coefficients of connection $\Pi^k_{ij}$ which are determined by the
     functions $P,Q,R$.
\\[1mm]
      In general case such type of connection is not metrizable  and corresponding space
      is not a Riemannian.

        Further we apply the notion
     of the Riemann extension
      of nonriemannian  space which were used earlier in  \cite{dryuma1:dryuma,dryuma2:dryuma,dryuma3:dryuma}.

     Remind basic properties of this construction.

     With help of the coefficients of affine connection of a given n-dimensional space
      can be introduced a 2n-dimensional
     Riemann space $D^{2n}$ with local coordinates $(x^i,\Psi_i)$ having the metric in form

\begin{equation} \label{dryuma:eq8}
{^{6}}ds^2=-2\Pi^k_{ij}(x^l)\Psi_k dx^i dx^j+2d \Psi_k dx^k
\end{equation}

\noindent where $\Psi_{k}$ are the additional coordinates.

The important property of such type metric is that the geodesic
 equations of metric (\ref{dryuma:eq8})  decomposes into two parts
\begin{equation} \label{dryuma:eq9}
\ddot x^k +\Gamma^k_{ij}\dot x^i \dot x^j=0,
\end{equation}
and
\begin{equation} \label{dryuma:eq10}
\frac{\delta^2 \Psi_k}{ds^2}+R^l_{kji}\dot x^j \dot x^i \Psi_l=0,
\end{equation}
where
\[
\frac{\delta \Psi_k}{ds}=\frac{d
\Psi_k}{ds}-\Pi^l_{jk}\Psi_l\frac{d x^j}{ds}
\]
and $R^l_{kji}$ are the curvature tensor of 3-dimensional space
with a given affine connection.

 The first part (\ref{dryuma:eq9})
is the system of equations for geodesic of basic space with local
coordinates $x^i$ and they does not contains the supplementary
coordinates $\Psi_k$.

 The second part (\ref{dryuma:eq10}) of the full system of geodesics  has the form
of linear $3\times3$ matrix system of second order ODE's for
supplementary  coordinates $\Psi_k$
\begin{equation} \label{dryuma:eq11}
\frac{d^2 \vec \Psi}{ds^2}+A(s)\frac{d \vec \Psi}{ds}+B(s)\vec \Psi=0.
\end{equation}

It is important to note that the geometry of extended space
connects with geometry of basic space. For example the property of
the space to be Ricci-flat keeps also for the extended space.

This fact give us the possibility to use the linear system of
equation (\ref{dryuma:eq11}) for the studying of geometrical
properties of the basic space.

In particular the invariants of  $k\times k$ matrix-function
\[
E=B-\frac{1}{2}\frac{d A}{ds}-\frac{1}{4}A^2
\]
under change of the coordinates $\Psi_k$ can be of used for that.

 The first applications the notion of extended spaces for the studying of nonlinear second order differential
 equations connected with nonlinear dynamical systems have been considered in the works of author
\cite{dryuma1:dryuma,dryuma2:dryuma,dryuma3:dryuma}.

    Here we consider some properties of the space defined by the Pfaff equations
     connected with a nonlinear dynamical systems.

\section{The  Lorenz dynamical system}

The equations of Lorenz dynamical system are

\begin{equation} \label{dryuma:eq12}
\frac{dx}{ds}=\sigma(y-x),\quad \frac{dy}{ds}=rx-y-zx,\quad
\frac{dz}{ds}=xy-bz,
\end{equation}
\medskip
where $\sigma,b,r$ are the parameters.

     These equations describe the behaviour of the flow lines along the vector field
\[
\vec N = [\sigma(y-x), rx-y-zx, xy-bz]
\]
depending from the parameters $\sigma,b,r$.

     The properties of  non holonomic variety $L^2$ in this case  are determined by the
     following Pfaff equation
\begin{equation} \label{dryuma:eq13}
\sigma(y-x)\frac{dx}{ds}+(rx-y-zx)\frac{dy}{ds}+(xy-bz)
\frac{dz}{ds}=0.
\end{equation}

   The object of holonomicity  for the Lorenz vector field is
\[
\left(\vec N, rot \vec
N\right)=\sigma\,xy-2\,\sigma\,{x}^{2}+{y}^{2}-bzr+b{z}^{2}+bz\sigma.
\]

     The properties of asymptotic lines of the
     corresponding variety $M^2$ are defined by the system of equations
$$
\left ({\frac {\partial }{\partial x}}P(x,y,z)\right )\left ({\frac {d
}{ds}}x(s)\right )^{2}\!+\!\left ({\frac {\partial }{\partial y}}Q(x,y,z) \right)
\left ({\frac {d}{ds}}y(s)\right )^{2}\!+\!\left ({\frac {
\partial }{\partial z}}R(x,y,z)\right )\times
$$
$$
\times\left ({\frac {d}{ds}}z(s)
\right )^{2}+\left ({\frac {\partial }{\partial
y}}P(x,y,z)+{\frac {\partial }{\partial x}}Q(x,y,z)\right )\left ({\frac {d}{ds}}x(s)
\right ){\frac {d}{ds}}y(s)+
$$
$$
+\left ({\frac {\partial }{\partial z}}P(x,
y,z)\!+\!{\frac {\partial }{\partial x}}R(x,y,z)\right )\left ({\frac {d}{ds}}x(s)
\right ){\frac {d}{ds}}z(s)+
$$
\begin{equation}
\label{dryuma:eq14}
+\left ({\frac {\partial }
{\partial y}}R(x,y,z)+{\frac {\partial }{\partial z}}Q(x,y,z)\right )
\left ({\frac {d}{ds}}z(s)\right ){\frac {d}{ds}}y(s) =0,
\end{equation}
\medskip

\begin{equation} \label{dryuma:eq15}
P(x,y,z){\frac {d}{ds}}x(s)+Q(x,y,z){\frac {d}{ds}}y(s)+R(x,y,z){
\frac {d}{ds}}z(s) =0,\end{equation}
\medskip
 which take the form
$$
-\left ({\frac {d}{ds}}y(s)\right )^{2}+\left (r+\sigma-z\right )
\left ({\frac {d}{ds}}x(s)\right ){\frac
{d}{ds}}y(s)-\sigma\,\left ({ \frac {d}{ds}}x(s)\right
)^{2}+
$$
\begin{equation}\label{dryuma:eq16}
+y\left ({\frac {d}{ds}}x(s)\right ){ \frac
{d}{ds}}z(s)-b\left ({\frac {d}{ds}}z(s)\right )^{2}=0,
\end{equation}
\medskip

\begin{equation}\label{dryuma:eq17}
\left (-bz+xy\right ){\frac {d}{ds}}z(s)+\left (rx-y-zx\right
){\frac {d}{ds}}y(s)+\left (\sigma\,y-\sigma\,x\right ){\frac
{d}{ds}}x(s) =0.
 \end{equation}
\\[1mm]

    We find from these equations the
$$
{\frac {d}{ds}}x(s)=-\displaystyle{\frac{1}{\sigma\,\left (-y+x\right )}}\times
$$
\begin{equation}\label{dryuma:eq18}
\times \left \{ -\left ({\frac {d}{ds}}y(s)\right)rx
+ \left ({\frac {d}{ds}}y(s)\right )y+\left ({\frac
{d}{ds}}y(s)\right ) zx-\left ({\frac {d}{ds}}z(s)\right )xy+
\left({\frac {d}{ds}}z(s) \right )bz \right \},
\end{equation}
and
\begin{equation}\label{dryuma:eq19}
 A\left (x(s)\right )^{2}+Bx(s)+C=0,
\end{equation}
\medskip
\noindent where

$$
A=\left (-\sigma\,z(s)+\sigma\,r-\sigma\right )\left ({\frac
{d}{ds}}y (s)\right )^{2}+\left ({\frac {d}{ds}}y(s)\right
)\sigma\,\left ({ \frac {d}{ds}}z(s)\right )y(s)-b\left ({\frac
{d}{ds}}z(s)\right )^{2} \sigma ,
$$

$$
B=\left (-\left (y(s)\right )^{3}\!+\!2\,b\sigma\,y(s)\!+\!y(s)bz(s)\right ) \left
({\frac {d}{ds}}z(s)\right )^{2}+
$$
$$
+\left (rbz(s)\!+\!2\,z(s)\left (y(s )\right)^{2}\!-\!\sigma\,bz(s)\!-\!
\left (z(s)\right )^{2}b\!-\!2\,r\left (y(s) \right)^{2}\!-\!\sigma\,
\left (y(s)\right )^{2}\!+\!\left (y(s)\right )^{2} \right )\times
$$
$$
\times\left({\frac {d}{ds}}y(s)\right ){\frac {d}{ds}}z(s)+\left (-
z(s)y(s)\!+\!\sigma\,z(s)y(s)\!-\!{r}^{2}y(s)\!-\!\sigma\,ry(s)\!+\!2\,rz(s)y(s)\!+\right.
$$
$$
\left.+\!\sigma\,y(s)\!+\!ry(s)\!-\!\left (z(s)\right )^{2}y(s)\right )\left ({\frac {d
}{ds}}y(s)\right )^{2},
$$

$$
C=\left (-b\sigma\,\left (y(s)\right )^{2}-{b}^{2}\left
(z(s)\right )^ {2}+\left (y(s)\right )^{2}bz(s)\right )\left
({\frac {d}{ds}}z(s) \right )^{2}+
$$
$$
+\left (-\left (z(s)\right)^{2}by(s)+rbz(s)y(s)-2\,y(s)bz (s)+\left (y(s)\right
)^{3}+\sigma\,bz(s)y(s)\right )\times
$$
$$
\times\left ({\frac {d} {ds}}y(s)\right){\frac {d}{ds}}z(s)+\left (r\left (y(s)\right )^{2}
- \left(y(s)\right )^{2}-z(s)\left (y(s)\right )^{2}\right )\left ({
\frac {d}{ds}}y(s)\right )^{2}.
$$
\medskip

\noindent After differentiating the equation (\ref{dryuma:eq19}) on the variable $s$ and
taking into account the expression (\ref{dryuma:eq18}) for
$\displaystyle\frac{dx(s)}{ds}$ we  get the relation

\begin{equation}\label{dryuma:eq20}
 E\left (x(s)\right )^{3}+F\left (x(s)\right )^{2}+Hx(s)+K=0
\end{equation}
\medskip

\noindent with some functions $E,F,H,K$ which does not contains
the variable $x(s)$.

     The resultant of the equations (\ref{dryuma:eq19}) and (\ref{dryuma:eq20})
     with regard of the variable
$$
\displaystyle{\frac {d}{ds}}z(s)=
\displaystyle{\frac {\displaystyle{\frac {d}{ds}}y(s)}{\displaystyle{\frac{d}{dz}}y(z)}}
$$
give us a following conditions

\begin{equation}\label{dryuma:eq21}
\left (r-1-z\right )y(z){\frac {d}{dz}}y(z)-bz+\left (y(z) \right
)^{2} =0,
\end{equation}
\smallskip

\begin{equation}\label{dryuma:eq22}
L\left ({\frac {d}{dz}}y(z)\right )^{2}+M{\frac {d}{dz}}y(z)+N=0,
\end{equation}
where  the coefficients of equation are
\[
{\frac {L}{1-r+z}}=-{z}^{4}{b}^{2}+\left
(2\,\sigma\,{b}^{2}+b\left (y (z)\right )^{2}+2\,r{b}^{2}\right
){z}^{3}+\]\[+\left (\left (-2\,rb-2\, \sigma\,b\right )\left
(y(z)\right )^{2}+4\,\sigma\,{b}^{2}-{r}^{2}{b}
^{2}-{\sigma}^{2}{b}^{2}-2\,\sigma\,{b}^{2}r\right
){z}^{2}+\]\[+\left ( \left (y(z)\right )^{4}+\left
({\sigma}^{2}b+2\,\sigma\,rb+{r}^{2}b-4 \,\sigma\,b\right )\left
(y(z)\right )^{2}\right )z+\left (1-r\right ) \left (y(z)\right
)^{4} ,
\]
\\[1mm]
\[
-1/2\,{\frac {M}{y(z)\left (1-r+z\right
)}}=-2\,{z}^{3}{b}^{2}+\left (
3\,\sigma\,{b}^{2}+3\,r{b}^{2}+b\left (y(z)\right )^{2}\right
){z}^{2} +\]\[+\left (\left (-2\,\sigma\,b-rb-2\,b\right )\left
(y(z)\right )^{2}-2
\,\sigma\,{b}^{2}r-{\sigma}^{2}{b}^{2}+4\,\sigma\,{b}^{2}-{r}^{2}{b}^{
2}\right )z+\]\[+\left (y(z)\right )^{4}+\left
(-3\,\sigma\,b+\sigma\,rb+rb +{\sigma}^{2}b\right )\left
(y(z)\right )^{2} ,
\]
\\[1mm]
\[
N={b}^{3}{z}^{4}+\left (-3\,\left (y(z)\right
)^{2}{b}^{2}-2\,r{b}^{3} +2\,\sigma\,{b}^{3}\right
){z}^{3}+\]\[+\left (b\left (y(z)\right )^{4}+ \left
(5\,r{b}^{2}+4\,\sigma\,{b}^{2}\right )\left (y(z)\right )^{2}+{
\sigma}^{2}{b}^{3}+{r}^{2}{b}^{3}-2\,\sigma\,{b}^{3}r\right
){z}^{2}+\]\[+ \left (\left (-rb-2\,\sigma\,b-3\,b\right )\left
(y(z)\right )^{4}+ \left
(-8\,\sigma\,{b}^{2}r-2\,{r}^{2}{b}^{2}+12\,\sigma\,{b}^{2}-2\,{
\sigma}^{2}{b}^{2}\right )\left (y(z)\right )^{2}\right
)z+\]\[+\left (y(z) \right )^{6}+\left
(2\,rb+{\sigma}^{2}b+2\,\sigma\,rb-b-4\,\sigma\,b \right )\left
(y(z)\right )^{4}+\left (4\,\sigma\,{b}^{2}+4\,\sigma\,{
r}^{2}{b}^{2}-8\,\sigma\,{b}^{2}r\right )\left (y(z)\right )^{2}
\]

    The solutions of the first order differential equations
    (\ref{dryuma:eq21}), (\ref{dryuma:eq22})
    together with conditions (\ref{dryuma:eq16})-(\ref{dryuma:eq17}) allow us to get the expressions
    for asymptotic line of the variety  $L^2$.

     Let us consider some examples.

     The  equation (\ref{dryuma:eq22}) has a set of singular
     solutions $y(z)$.

One of them determines from the relation
\[
{z}^{2}b+\left (-2\,\sigma\,b-2\,rb\right
)z+{y}^{2}+2\,br\sigma+b{r}^{2}-4\, \sigma\,b+b{\sigma}^{2}=0
\]
which presents the second order curve in the plane (z,y).

   Remark that the equation (\ref{dryuma:eq22}) presents the algebraic
   curve
\[\Phi(y',y,z)=0\]
 of genus $g=1$ with respect to variables  $y',y$ in case $r\neq 1$ and genus $g=0$ when $r=1$.

 According with general theory some of such type equations can be
 integrated with help of elliptic functions or can be brought to integration of the Rikkati equation.

     In both cases the  properties  of
     asymptotic lines should be dependent from the parameters of model.

\section{The  R\"ossler dynamical system}

   Differential  equations of the R\"ossler dynamical model are

\begin{equation} \label{dryuma:eq23}
\frac{dx}{ds}=-y-z,\quad \frac{dy}{ds}=x+ay,\quad
\frac{dz}{ds}=b+xz-cz,
\end{equation}
where $a,b,c$ are the parameters.

     For corresponding non holonomic variety $V^2$ we have a following Pfaff equation
\begin{equation} \label{dryuma:eq24}
(-y-z)\frac{dx}{ds}+(x+ay)\frac{dy}{ds}+(b+xz-cz) \frac{dz}{ds}=0.
\end{equation}

   The object of holonomicity  for the R\"ossler vector field is
\[
\left(\vec N, rot \vec N\right)=-x+xz-ay-ayz+2\,b-2\,cz.
\]

     In this case the properties of the system (\ref{dryuma:eq16})-(\ref{dryuma:eq17})
     for asymptotic lines are determined by the equation
\begin{equation} \label{dryuma:eq25}
A{\frac {d^{2}}{d{z}^{2}}}y(z)+B =0
\end{equation}
where
$$
A=\left (-y(z)az+az-a{z}^{2}+y(z)a\right )\left ({\frac{d}{dz}}y(z) \right )^{2}+
$$
$$
+\left(-2\,a{z}^{3}-2\,y(z)az-2\,a{z}^{2}y(z)-2\,\left (y (z)\right)^{2}a\right )
{\frac {d}{dz}}y(z)-
$$
$$
-ay(z){z}^{3}-2\,bz-{z}^{3}
c+c{z}^{2}+y(z)c+b{z}^{2}-cy(z)z-a\left (y(z)\right )^{2}z+
$$
$$
+\left(y(z) \right )^{2}a+b+a{z}^{2}y(z)
$$
and
$$
B=2\left (z+y(z)\right )\left (\left ({\frac {d}{dz}}y(z)\right)
ay(z) +\left ({\frac {d}{dz}}y(z)\right )c-\right.
$$
$$
\left.-\left ({\frac{d}{dz}}y(z)\right )^{3}a-az\left ({\frac {d}{dz}}y(z)\right)^{2}+b\right ).
$$

   This is equation of the second
   range at the condition $a \neq 0$.

   In the case $a=0$ it takes a form of the first range equation
$$
\left( 2\,bz\!-\!y \left( z \right) c\!+\!{z}^{3}c\!-\!c{z}^{2}\!-\!b\!+\!cy \left( z
\right) z\!-\!b{z}^{2} \right) {\frac {d^{2}}{d{z}^{2}}}y \left( z
\right)\!+\!2\,cz{\frac {d}{dz}}y \left( z \right)\!+
$$
$$
+\!2\,y \left(z\right) c{\frac {d}{dz}}y \left(z\right)\!+\!2\,bz+2\,y \left(z\right) b=0.
$$
    At the conditions $a=0,\quad b=0$, $c \neq0$ the properties of
    asymptotic lines of corresponding non holonomic variety are dependent
     from the solutions of the equation
    \[
 \left( {z}^{3}-{z}^{2}-y \left( z \right) +y \left( z \right) z
 \right) {\frac {d^{2}}{d{z}^{2}}}y \left( z \right) +2\,z{\frac {d}{d
z}}y \left( z \right) +2\,y \left( z \right) {\frac {d}{dz}}y
\left( z
 \right)=0.
\]

\section{Quadratic systems}

    Here we consider the properties of the Pfaff equations
    connected with a polynomial differential systems in $R^2$
    defined by
\begin{equation} \label{dryuma:eq26}
\frac{dx}{ds}=p(x,y),\quad \frac{dy}{ds}=q(x,y),
\end{equation}
where $p(x,y)$ and $q(x,y)$ are polinomials of degree 2.

    The system (\ref{dryuma:eq26}) takes the form of the Pfaff equation
    after extension on a projective plane
\begin{equation}\label{dryuma:eq27}
 \left( x\,Q \left( x,y,z \right) -y\,P \left( x,y,z \right) \right) {
\it dz}-z\,Q \left( x,y,z \right){\it dx} +z\,P \left( x,y,z
\right) { \it dy} =0,
\end{equation}
where
\[
P(x,y,z),\quad Q(x,y,z),\quad R(x,y,z)
\]
are the homogeneous functions constructed from the functions
$p(x,y)$ and $q(x,y)$.

       As example for the system
\[
\frac{dx}{ds}=kx+ly+a{x}^{2}+bxy+c{y}^{2} ,\quad \frac{dy}{ds}=
=mx+ny+e{x}^{2}+fxy+h{y}^{2},
\]
with a ten parameters one get a Pfaff equation after a projective extension
\[
P(x,y,z)dx+ Q(x,y,z)dy+R(x,y,z)dz=0,
\]
or
 $$
 \left(  \left( -ny-mx \right) {z}^{2}+ \left( -h{y}^{2}-fxy-e{x}^{2}
 \right) z \right) {\it dx}+ \left( \left( ly+kx \right) {z}^{2}+\right.
 $$
 $$
\left.+ \left( c{y}^{2}+bxy+a{x}^{2} \right) z \right) {\it dy}+ \left(
\left( -l{y}^{2}+ \left( n-k \right) xy+m{x}^{2} \right) z-c{y}^{3}+\right.
$$
\begin{equation} \label{dryuma:eq28}
\left.+\left( h-b \right) x{y}^{2}+ \left( f-a \right) {x}^{2}y+e{x}^{3}
\right) {\it dz}=0.
\end{equation}

    This equation  corresponds the system
\begin{equation} \label{dryuma:eq29}
\frac{dx}{ds}=P(x,y,z),\quad \frac{dy}{ds}=Q(x,y,z),\quad\frac{dz}{ds}=R(x,y,z),
\end{equation}
 where
\[
P(x,y,z)=  \left( -ny-mx \right) {z}^{2}+ \left( -h{y}^{2}-fxy-e{x}^2\right) z  ,
\]
\[
Q(x,y,z)=  \left( ly+kx \right) {z}^{2}+
 \left( c{y}^{2}+bxy+a{x}^{2} \right) z,
 \]
\[
R(x,y,z)=  \left(-l{y}^{2}+ \left( n-k \right) xy+m{x}^{2} \right) z-c{y}^{3}+
 \left( h-b \right) x{y}^{2}+ \left( f-a \right) {x}^{2}y+e{x}^{3}
\]

     It is important to note that the vector field $\vec N=(P,Q,R)$ connected with
     a system (\ref{dryuma:eq29}) is holonomic due the condition
\[
(\vec N,rot \vec N)=0.
\]

     Corresponding Pfaff (\ref{dryuma:eq28}) equation is integrable and determines
     the family of developing surfaces
$U(x,y,z)=C$ in the $R^3$-space.

      The investigation of the asymptotic lines of
      the surfaces which corresponds  the equation
      (\ref{dryuma:eq28}) may be useful in applications.

  Let us consider some of examples.

  The system of equations
\[
\frac{dx}{ds}=5x+6{x}^{2}+4(1+\mu)xy+\mu{y}^{2} ,\quad
\frac{dy}{ds}= =x+2y+4xy+(2+3\mu){y}^{2},
\]
with $(-71+17\sqrt(17)/32< \mu < 0$ posseses the invariant
algebraic curve
\[
x^2+x^3+x^2y+2\mu xy^2+2\mu xy^3+mu^2 y^4=0.
\]

       The conditions of compatibility of equations for asymptotic lines
of the variety defined by the corresponding  Pfaff equation lead
to a following  conditions on the functions

\[
 \left( 3\,{\mu}^{2}+4\,\mu+2 \right)  \left( y \left( s \right)
 \right) ^{2}+ \left( 9\,\mu+6 \right) x \left( s \right) y \left( s
 \right) +6\, \left( x \left( s \right)  \right) ^{2}=0,
\]
\[
 \left( y \left( s \right)  \right) ^{2}\mu+ \left( -2\,\mu-1 \right)
x \left( s \right) y \left( s \right) -3\, \left( x \left( s
\right)
 \right) ^{2}
=0.
\]

   The values of parameter $\mu$ determined by the conditions
\[
(3\mu+10)(\mu+4)(\mu-2)=0,\quad \mu=0
\]
are special.

    Next example is the system with at least a four limit cycles.
\[
\frac{dx}{ds}=kx-y-10x^2+bxy+y^2 ,\quad \frac{dy}{ds}= =x+x^2+fxy,
\]

   The studying of asymptotic lines for this system give a following
   conditions on parameters
\[
10\,{k}^{2}+11\,f{k}^{2}+{f}^{2}{k}^{2}-7\,k+2\,fk-9\,b+{k}^{3}-80-81
\,f+{b}^{2}k+2\,b{k}^{2}+9\,bk+9\,fbk=0
\]
and functions
\[{x}^{3}+10\,y \left( x \right) {x}^{2}- \left( y \left( x \right)
 \right) ^{3}+f{x}^{2}y \left( x \right) -bx \left( y \left( x
 \right)  \right) ^{2}=0,
\]
\[f \left( y \left( x \right)  \right) ^{2}-fy \left( x \right) kx+
 \left( y \left( x \right)  \right) ^{2}+y \left( x \right) bx+y
 \left( x \right) x-k{x}^{2}-10\,{x}^{2}=0.
\]

    Remark that these equations  equivalent the equations of direct
    lines.

\section{Cubic systems}

    By analogy can be investigated the properties of the asymptotic lines
    of the surfaces  connected with the system
\begin{equation} \label{dryuma:eq27}
\frac{dx}{ds}=p(x,y),\quad \frac{dy}{ds}=q(x,y),
\end{equation}
where $p(x,y)$ and $q(x,y)$ are polinomials of degree 3.

    Let us consider the examples.

 The system
\begin{equation} \label{dryuma:eq28}
\frac{dx}{ds}=y,\quad \frac{dy}{ds}=-x-x^2y+mu^2y
\end{equation}
connects with the Van der Pol equation.

After extension on the projective plane we get an integrable Pfaff
equation
\[
 \left( -{x}^{2}{z}^{2}-{x}^{3}y+x{\mu}^{2}y{z}^{2}-{y}^{2}{z}^{2}
 \right) {\it dz}+ \left( {z}^{3}x+z{x}^{2}y-{z}^{3}{\mu}^{2}y
 \right) {\it dx}+y{z}^{3}{\it dy}=0.
\]

   The equations for the asymptotic lines of corresponding surface
   give the conditions
\[{x}^{2}+ \left( y \left( x \right)  \right) ^{2}-x{\mu}^{2}y \left( x
 \right)=0
\]
and
\[
-{x}^{2}+2\, \left( y \left( x \right)  \right)
^{2}+4\,x{\mu}^{2}y
 \left( x \right)=0.
\]

    From the first condition we find
\[y(x)= \left( 1/2\,{\mu}^{2}+(-)1/2\,\sqrt {{\mu}^{4}-4} \right) x,
\]
and the second gives us
\[y(x)= \left( -{\mu}^{2}+(-)1/2\,\sqrt {4\,{\mu}^{4}+2} \right) x.
\]

 For the system
\begin{equation} \label{dryuma:eq29}
\frac{dx}{ds}=-y+ax(x^2+y^2-1),\quad \frac{dy}{ds}=x+by(x^2+y^2-1)
\end{equation}
the point $(0,0$ is node at the condition $ab>-1$ and $(a-b)^2>=4$
and it has the limit cycle around this point.

After extension we get an integrable Pfaff equation
$$
\left( {x}^{2}{z}^{2}+by{x}^{3}+xb{y}^{3}-xby{z}^{2}+{y}^{2}{z}^{2}-y
a{x}^{3}-ax{y}^{3}+yax{z}^{2} \right) {\it dz}+
$$
$$
+\left(-{z}^{3}y+za{x} ^{3}+zax{y}^{2}-{z}^{3}ax \right) {\it dy}+
\left( -{z}^{3}x-zby{x}^{2 }-zb{y}^{3}+{z}^{3}by \right) {\it dx}=0
$$
with a developing surfaces as general integral.

 The conditions on parameters for existence of  asymptotic lines are
$$
ab+1=0,\qquad -a+b=0.
$$

   For the corresponding functions one get the equations of direct lines

\[xby \left( x \right) -{x}^{2}-y \left( x \right) ax- \left( y \left( x
 \right)  \right) ^{2}=0,
\]
or
\[y(x)=xby \left( x \right) -{x}^{2}-y \left( x \right) ax- \left( y \left( x
 \right)  \right) ^{2},
\]

and
$$
\left( 30\,ab+8\,{a}^{3}b-25\,{a}^{2}-9\,{b}^{2}-12\,{a}^{4} \right) {x}^{6}+
$$
$$
+\left( 36\,{b}^{3}+76\,b{a}^{2}-108\,a{b}^{2}-4\,{a}^{3}-32\,
{a}^{3}{b}^{2}+32\,b{a}^{4} \right) y{x}^{5}+
$$
$$
+\left(69\,{b}^{2}+53\,{a}^{2}+36\,{a}^{4}-64\,{a}^{3}b-134\,ab+24\,a{b}^{3}-8\,{a}^{2}{b}^{2}
 \right) {y}^{2}{x}^{4}+
 $$
 $$
 +\left( -24\,{b}^{3}-8\,b{a}^{2}+64\,{a}^{3}{b}^{2}+8\,a{b}^{2}+24\,{a}^{3}-64\,{a}^{2}{b}^{3} \right)
{y}^{3}{x}^{3 }+
$$
$$
+\left(-8\,{a}^{2}{b}^{2}-134\,ab+24\,{a}^{3}b+36\,{b}^{4}-64\,a{b}^{3}+69\,
{a}^{2}+53\,{b}^{2} \right) {y}^{4}{x}^{2}+
$$
$$
+ \left(-32\,{b}^{4}a-36\,{a}^{3}-76\,a{b}^{2}+32\,{a}^{2}{b}^{3}+4\,{b}^{3}+108\,
b{a}^{2} \right) {y}^{5}x+
$$
$$
+\left(-25\,{b}^{2}-9\,{a}^{2}+30\,ab+8\,a{b}^{3} -12\,{b}^{4} \right){y}^{6}=0
$$
with a very complicate relations between the parameters $a,b$.

\section{Quatric systems}

    The system
\begin{equation} \label{dryuma:eq30}
\frac{dx}{ds}=-Ay+yx^2-x^4,\quad \frac{dy}{ds}=ax-x^3
\end{equation}
has a center and the limit cycle at the conditions $a=2A+A^2$,
 $0 < A < 5.10^{-5}$.

After extension we get an integrable Pfaff equation
$$
\left( {z}^{2}{x}^{3}-{z}^{4}ax \right) {\frac {d}{ds}}x \left( s \right)+
\left( -{x}^{4}z-{z}^{4}Ay+{z}^{2}{x}^{2}y \right) {\frac {d}{ds}}y \left( s \right) +
$$
$$
+\left(a{x}^{2}{z}^{3}-{x}^{4}z+A{y}^{2}{z }^{3}-z{x}^{2}{y}^{2}+y{x}^{4}
\right) {\frac {d}{ds}}z \left( s \right) =0
$$
with a developing surfaces as  general integral.

   For the functions one get the equations of direct lines

\[{x}^{10}+ \left( 3\,a-9\,A \right) {z}^{2}{x}^{8}+ \left( 24\,{A}^{2}-
3\,Aa \right) {z}^{4}{x}^{6}+ \left(
{A}^{2}-36\,{A}^{3}-2\,Aa+{a}^{2}
 \right) {z}^{6}{x}^{4}+ \]\[\left( 6\,{A}^{3}+36\,{A}^{3}a-12\,{A}^{2}a+6
\,{a}^{2}A \right) {z}^{8}{x}^{2}+ \left(
9\,{a}^{2}{A}^{2}-18\,{A}^{3 }a+9\,{A}^{4} \right) {z}^{10} =0
\]
and the family of conics
\[ \left( z \left( s \right)  \right) ^{2}a- \left( x \left( s \right)
 \right) ^{2}=0
\]
\[
 \left( x \left( s \right)  \right) ^{2}-A \left( z \left( s \right)
 \right) ^{2}
=0.
\]

\section{Geodesics of the first kind}

   A geodesics of the first kind on non holonomic variety are defined by the system of equations
(\ref{dryuma:eq5})-(\ref{dryuma:eq7}) .

   Let us consider the solutions of this system of equations for the
   variety $V^2$ defined by the vector field
\[
    \vec F=[ayz,bxz,cxy],
\]
where $a,b,c$ are parameters.

    In this case the condition  $\vec N=0$ holds and variety is holonomic.

  The system of equations for geodesics is
$$
\displaystyle{\displaystyle\frac {d^{2}}{d{s}^{2}}}x(s)+\displaystyle{\frac {ayzx\left (c+b\right )\left(
\displaystyle{ \frac {d}{ds}}y(s)\right )\displaystyle{\frac{d}{ds}}z(s)}{{a}^{2}{y}^{2}{z}^{2}+
{b}^{2}{z}^{2}{x}^{2}+{c}^{2}{x}^{2}{y}^{2}}}+
\displaystyle{\frac{a{y}^{2}z\left (a +c\right )\left (\displaystyle{\frac {d}{ds}}x(s)\right)
\displaystyle{\frac {d}{ds}}z(s)}{{a}^{2}{y}^{2}{z}^{2}+{b}^{2}{z}^{2}{x}^{2}+{c}^{2}{x}^{2}{y}^{2}}}+
$$
$$
+\displaystyle{\frac{ay{z}^{2}\left (a+b\right )\left (\displaystyle{\frac {d}{ds}}x(s)\right)
\displaystyle{\frac {d}{ds}}y(s)}{{a}^{2}{y}^{2}{z}^{2}+{b}^{2}{z}^{2}{x}^{2}+{c}^{2}{x}^{2}{y}^{2}}} =0,
$$

$$
\displaystyle{\frac {d^{2}}{d{s}^{2}}}y(s)+\displaystyle{\frac {b{x}^{2}z\left (c+b\right )
\left ({\frac {d}{ds}}y(s)\right )\displaystyle{\frac{d}{ds}}z(s)}{{a}^{2}{y}^{2}
{z}^{2}+{b}^{2}{z}^{2}{x}^{2}+{c}^{2}{x}^{2}{y}^{2}}}+
\displaystyle{\frac {bxzy\left (a+c\right )\left (\displaystyle{\frac {d}{ds}}x(s)\right )
\displaystyle{\frac{d}{ds}}z(s)}{{a}^{2}{y}^{2}{z}^{2}+{b}^{2}{z}^{2}{x}^{2}+{c}^{2}{x}^{2}{y}^{2}}}+
$$
$$
+{\frac {bx{z}^{2}\left (a+b\right )\left ({\frac{d}{ds}}x(s)\right )\displaystyle {\frac
{d}{ds}}y(s)}{{a}^{2}{y}^{2}{z}^{2}+{b}^{2}{z}^{2}{x}^{2}+{c}^{2}{x}^{2}{y}^{2}}} =0,
$$

$$
\displaystyle{\frac {d^{2}}{d{s}^{2}}}z(s)+\displaystyle{\frac {c{x}^{2}y\left (c+b\right )
\left (\displaystyle{\frac {d}{ds}}y(s)\right )\displaystyle{\frac{d}{ds}}z(s)}{{a}^{2}{y}^{2}
{z}^{2}+{b}^{2}{z}^{2}{x}^{2}+{c}^{2}{x}^{2}{y}^{2}}}+\displaystyle{\frac{cx{y}^{2}
\left (a+c\right )\left (\displaystyle{\frac {d}{ds}}x(s)\right)\displaystyle{\frac {d}{ds}}
z(s)}{{a}^{2}{y}^{2}{z}^{2}+{b}^{2}{z}^{2}{x}^{2}+{c}^{2}{x}^{2}{y}^{2}}}+
$$
$$
+\displaystyle{\frac {cxyz\left (a+b\right )\left (\displaystyle{\frac{d}{ds}}x(s)\right )
\displaystyle{ \frac{d}{ds}}y(s)}{{a}^{2}{y}^{2}{z}^{2}+{b}^{2}{z}^{2}{x}^{2}+{c}^{2}{x}^{2}{y}^{2}}} =0.
$$

    This system of equations has an integral
\[
ay(s)z(s){\frac {d}{ds}}x(s)+bx(s)z(s){\frac
{d}{ds}}y(s)+cx(s)y(s){ \frac {d}{ds}}z(s) =0
\]
and can be integrated.

    In fact after substitution of the expression ${ \frac
    {d}{ds}}z(s)$ into the above equations one get the algebraic relations
    with respect the variable $z(s)$ which are compatible at the
    conditions
\[
\left (ab+{b}^{2}\right )\left (z(y)\right )^{2}+2\,bz(y)cy{\frac
{d}{ dy}}z(y)+\left (a{y}^{2}c+{c}^{2}{y}^{2}\right )\left ({\frac
{d}{dy}} z(y)\right )^{2} =0
\]
and lead to the solution
\[
z(y)={y}^K{\it \_C1}
\]
where
 \[
K={-{\frac {cb}{{c}^{2}+ac}}-{\frac {\sqrt {-cab\left (a+c+b
\right )}}{{c}^{2}+ac}}}.
\]

   Finally  we present the expression for the Chern-Simons invariant of affine
   connection defined by the equations of geodesics (\ref{dryuma:eq5}-\ref{dryuma:eq7})
\begin{equation} \label{dryuma:eq34}
CS=\int\epsilon^{i j k}(\Gamma^p_{i q}\Gamma^q_{k
p;j}+\frac{2}{3}\Gamma^p_{i q}\Gamma^q_{j r}\Gamma^r_{k p})dx dy
dz.
\end{equation}

   In the case of  the Lorenz system of equations it has been obtained with the help of
   a six dimensional Riemann extension of  corresponding space and has the form
\[
CS(\Gamma)=\int\frac{L}{2M^2}dxdydz,
\]
where
\[
L=  \left( 2\,b+2-2\,\sigma \right) {x}^{2}{y}^{2}+ \left( 3\,{\sigma}^{2
}+4\,\sigma\,r-4\,rb-2\,b\sigma-4\,r \right) z{x}^{2}+\left( 2\,b+2-2 \,\sigma
\right) {z}^{2}{x}^{2}+ \]\[+\left( -3\,r{\sigma}^{2}+4\,{\sigma}^
{2}b-5\,{\sigma}^{3}+2\,b\sigma\,r-2\,\sigma\,{r}^{2}\!+\!4\,{\sigma}^{2}+
2\,b{r}^{2}+2\,{r}^{2} \right) {x}^{2}+\]\[+ \left( -2\,rb-4\,r+9\,{\sigma}
^{3}+2\,r{\sigma}^{2} \right) yx+ \left( -2\,\sigma+4-2\,{\sigma}^{2}-
2\,\sigma\,r+2\,b\sigma-4\,{b}^{2} \right) zyx-\]\[- \left(  \left( b\sigma
\,r+2\,{\sigma}^{2}+{\sigma}^{2}b-\sigma\,{z}^{2}-2\,\sigma\,r-\sigma \,{r}^{2}
\right) y-\sigma\,{y}^{3} \right) x+ \left( -2\,\sigma\,{b}^
{2}+2\,{b}^{3}-2\,{b}^{2} \right) {z}^{2}+\]\[+ \left( -4\,{\sigma}^{3}-{
\sigma}^{2}b+2-2\,b-{\sigma}^{2}+\sigma\,r+ \left( -\sigma+b\sigma
 \right) z-b\sigma\,r-2\,\sigma \right) {y}^{2}+\]\[+ \left( -2\,\sigma\,{b
}^{2}+2\,r{b}^{2}+2\,{b}^{2}r\sigma-2\,{b}^{2}{\sigma}^{2} \right) z,
 \]
and
\[
M= \left( {x}^{2}+{b}^{2} \right) {z}^{2}+ \left(  \left( -2\,b+2
 \right) xy-2\,r{x}^{2} \right) z+\]\[+ \left( {\sigma}^{2}+1+{x}^{2}
 \right) {y}^{2}+ \left( -2\,{\sigma}^{2}-2\,r \right) xy+ \left( {
\sigma}^{2}+{r}^{2} \right) {x}^{2}.
\]


\begin{thebibliography}{10}

\bibitem{dryuma:sintsov} {\sc Sintsov D.M.} {\it Raboty po negolonomnoi geometrii}. Kiev, 1972.

\bibitem{dryuma4:sluchaev} {\sc Sluchaev V.V.} {\it Geometriya vektornyh polei}.
Tomsk, Izd. Tomskogo universiteta, 1982.

\bibitem{dryuma5:aminov} {\sc Aminov Yu.A.} {\it Geometriya vektornogo polya}. Moskva, Nauka, 1990.

\bibitem{dryuma1:dryuma} {\sc Dryuma V.} {\it The Riemann Extensions in theory of differential equations and their
applications}. Matematicheskaya fizika, analiz, geometriya, 2003, {\bf 10}, N~3, p.~1--19.

\bibitem{dryuma2:dryuma} {\sc Dryuma V.} {\it The Riemann and Einstein-Weyl geometries in the theory of ODE's,
their applications and all that}. New Trends in Integrability and Partial Solvability, p.~115-156
(eds. A.B.~Shabat et al.).  Kluwer Academic Publishers (ArXiv: nlin: SI/0303023, 11 March, 2003, p.~1--37).

\bibitem{dryuma3:dryuma} {\sc Dryuma V.} {\it Applications of Riemannian and Einstein-Weyl Geometry in
the theory of second order ordinary differential equations}. Theoretical and  Mathematical Physics,
2001, {\bf 128}, N~1, p.~845--855.

\end{thebibliography}
\end{document}